\documentclass[onecolumn]{article}
\usepackage{lmodern}
\usepackage{amssymb,amsmath}
\usepackage{ifxetex,ifluatex}
\usepackage{fixltx2e} 
\ifnum 0\ifxetex 1\fi\ifluatex 1\fi=0 
  \usepackage[T1]{fontenc}
  \usepackage[utf8]{inputenc}
\else 
  \ifxetex
    \usepackage{mathspec}
  \else
    \usepackage{fontspec}
  \fi
  \defaultfontfeatures{Ligatures=TeX,Scale=MatchLowercase}
\fi
\IfFileExists{upquote.sty}{\usepackage{upquote}}{}
\IfFileExists{microtype.sty}{%
\usepackage[]{microtype}
\UseMicrotypeSet[protrusion]{basicmath} 
}{}
\PassOptionsToPackage{hyphens}{url} 
\usepackage[unicode=true]{hyperref}
\hypersetup{
            pdftitle={Nonlinear Optimization in R using nlopt},
            pdfauthor={Rahul Bhadani},
            pdfborder={0 0 0},
            breaklinks=true}
\usepackage[margin=1in]{geometry}
\usepackage{natbib}
\usepackage{color}
\usepackage{fancyvrb}

\DefineVerbatimEnvironment{Highlighting}{Verbatim}{commandchars=\\\{\}}
\usepackage{framed}
\definecolor{shadecolor}{RGB}{248,248,248}
\newenvironment{Shaded}{\begin{snugshade}}{\end{snugshade}}
\newcommand{\KeywordTok}[1]{\textcolor[rgb]{0.13,0.29,0.53}{\textbf{#1}}}
\newcommand{\DataTypeTok}[1]{\textcolor[rgb]{0.13,0.29,0.53}{#1}}
\newcommand{\DecValTok}[1]{\textcolor[rgb]{0.00,0.00,0.81}{#1}}

\newcommand{\FloatTok}[1]{\textcolor[rgb]{0.00,0.00,0.81}{#1}}

\newcommand{\StringTok}[1]{\textcolor[rgb]{0.31,0.60,0.02}{#1}}

\newcommand{\CommentTok}[1]{\textcolor[rgb]{0.56,0.35,0.01}{\textit{#1}}}

\newcommand{\OtherTok}[1]{\textcolor[rgb]{0.56,0.35,0.01}{#1}}

\newcommand{\ControlFlowTok}[1]{\textcolor[rgb]{0.13,0.29,0.53}{\textbf{#1}}}
\newcommand{\OperatorTok}[1]{\textcolor[rgb]{0.81,0.36,0.00}{\textbf{#1}}}

\newcommand{\NormalTok}[1]{#1}
\usepackage{graphicx,grffile}
\makeatletter
\def\maxwidth{\ifdim\Gin@nat@width>\linewidth\linewidth\else\Gin@nat@width\fi}
\def\maxheight{\ifdim\Gin@nat@height>\textheight\textheight\else\Gin@nat@height\fi}
\makeatother
\setkeys{Gin}{width=\maxwidth,height=\maxheight,keepaspectratio}
\IfFileExists{parskip.sty}{%
\usepackage{parskip}
}{
\setlength{\parindent}{0pt}
\setlength{\parskip}{6pt plus 2pt minus 1pt}
}
\ifx\paragraph\undefined\else
\let\oldparagraph\paragraph
\renewcommand{\paragraph}[1]{\oldparagraph{#1}\mbox{}}
\fi
\ifx\subparagraph\undefined\else
\let\oldsubparagraph\subparagraph
\renewcommand{\subparagraph}[1]{\oldsubparagraph{#1}\mbox{}}
\fi

\makeatletter
\def\fps@figure{htbp}
\makeatother

\usepackage[none]{hyphenat}
\graphicspath{ {images/} }
\usepackage{color}
\definecolor{ocre}{RGB}{243,102,25}
\usepackage[font={color=ocre,bf},labelfont=bf]{caption}
\usepackage{hyperref}
\hypersetup{colorlinks=true,citecolor=blue}
\def\BibTeX{{\rm B\kern-.05em{\sc i\kern-.025em b}\kern-.08em T\kern-.1667em\lower.7ex\hbox{E}\kern-.125emX}}

\title{Nonlinear Optimization in R using \texttt{nlopt}}
\author{Rahul Bhadani\footnote{The University of Arizona,
  \href{mailto:rahulbhadani@email.arizona.edu}{\nolinkurl{rahulbhadani@email.arizona.edu}}}}
\date{10 January 2021}

\begin{document}
\maketitle
\begin{abstract}
In this article, we present a problem of nonlinear constraint
optimization with equality and inequality constraints. Objective
functions are defined to be nonlinear and optimizers may have a lower
and upper bound. We solve the optimization problem using the open-source
R package \texttt{nloptr}. Several examples have been presented.
\end{abstract}

\section{Introduction}\label{introduction}

Often in physical science research, we end up with a hard problem of
optimizing a function (called objective) that satisfies a range of
constraints - linear or non-linear equalities and inequalities. The
optimizers usually also have to adhere to the upper and lower bound. We
recently worked on a similar problem in Quantum Information Science
(QIS) where we attempted to optimize a non-linear function based on a
few constraints dictated by rules of physics and mathematics. Interested
readers may find our work on Constellation Optimization for Phase-Shift
Keying Coherent States \citep{bhadani2020constellation} where we
optimized mutual information for Quadri-Phase-Shift Keying (QPSK) based
on a set of constraints.

While chasing the problem of non-linear optimization with a set of
constraints, we found out that not all optimization routines are created
equally. There several libraries available in different languages such
as python (scipy.optimize), Matlab (fmincon), C++ (robotim, nlopt), and
R (nloptr). While the list for optimization routine presented here is
not exhaustive, some of them are more reliable than others, some provide
faster execution than others and some have better documentation. Based
on several key factors, we find nloptr, implemented in the R language to
be most suitable for nonlinear optimization. nloptr uses nlopt
implemented in C++ as a backend. As a result, it provides the elegance
of the R language and the speed of C++. The optimization procedure is
performed quickly in a fraction of seconds even with a tolerance of the
order of 10e-15.

\section{Nonlinear Optimization
Problem}\label{nonlinear-optimization-problem}

A general nonlinear optimization problem usually have the form

\[
\min_{x \in \mathbb{R}^n} f(x)
\] such that \[
g(x) \leq 0
\] \[
h(x) = 0
\] \[
x_L \leq x \leq x_U
\]

where \(f\) is an objective function, \(g\) defines a set of inequality
constraints, \(h\) is a set of equality constraints. \(x_L\) and \(x_U\)
are lower and upper bounds respectively. In the literature, several
optimization algorithms have been presented. For example, MMA (Method of
moving asymptotes) \citep{svanberg1987method} supports arbitrary
nonlinear inequality constraints, (COBYLA) Constrained Optimization BY
Linear Approximation \citep{powell1994direct}, (ORIG\_DRIECT) DIRECT
algorithm \citep{finkel2003direct}. Optimization algorithms that also
support nonlinear equality constraints include ISRES (Improved
Stochastic Ranking Evolution Strategy) \citep{runarsson2000stochastic},
(AUGLAG) Augmented Lagrangian Algorithm \citep{conn1991globally}. A full
list of such methods can be found on the nlopt C++ reference page at
\url{https://nlopt.readthedocs.io/en/latest/NLopt_Reference/}.

In the rest of the article, We provide several examples of solving a
constraint optimization problem using R. We use R Studio that combines R
compiler and editor. R Studio also provides a knitr tool which is great
for writing documentation or articles with inline code which can also
generate a latex source code and a pdf file. Most of the example
presented here has been modified from test suites used to validate
functions in \texttt{nloptr} R package.

\section{Installation and loading the
library}\label{installation-and-loading-the-library}

Installation of \texttt{nloptr} in R is fairly straightforward.

\begin{Shaded}
\begin{Highlighting}[]
\KeywordTok{install.packages}\NormalTok{(}\StringTok{"nloptr"}\NormalTok{)}
\end{Highlighting}
\end{Shaded}

\begin{Shaded}
\begin{Highlighting}[]
\KeywordTok{library}\NormalTok{(}\StringTok{'nloptr'}\NormalTok{)}
\end{Highlighting}
\end{Shaded}

\section{Example 1: Optimization with explicit
gradient}\label{example-1-optimization-with-explicit-gradient}

In the first example, we will minimize the Rosenbrock Banana function

\[
f(x) = 100 (x_2 - x_1^2)^2 + (1-x_1)^2
\] whose gradient is given by

\[
\nabla f(x) = \begin{pmatrix} -400x_1 (x_2 - x_1^2) - 2(1-x_1) \\ 200(x_2 - x_1^2)  \end{pmatrix}
\]

However, not all the algorithms in \texttt{nlopt} require explicit
gradient as we will see in further examples. Let's define the objective
function and its gradient first:

\begin{Shaded}
\begin{Highlighting}[]
\NormalTok{eval_f <-}\StringTok{ }\ControlFlowTok{function}\NormalTok{(x)}
\NormalTok{\{}
    \KeywordTok{return}\NormalTok{ ( }\DecValTok{100} \OperatorTok{*}\StringTok{ }\NormalTok{(x[}\DecValTok{2}\NormalTok{] }\OperatorTok{-}\StringTok{ }\NormalTok{x[}\DecValTok{1}\NormalTok{] }\OperatorTok{*}\StringTok{ }\NormalTok{x[}\DecValTok{1}\NormalTok{])}\OperatorTok{^}\DecValTok{2} \OperatorTok{+}\StringTok{ }\NormalTok{(}\DecValTok{1} \OperatorTok{-}\StringTok{ }\NormalTok{x[}\DecValTok{1}\NormalTok{])}\OperatorTok{^}\DecValTok{2}\NormalTok{ )}
\NormalTok{\}}


\NormalTok{eval_grad_f <-}\StringTok{ }\ControlFlowTok{function}\NormalTok{(x) \{}
\KeywordTok{return}\NormalTok{( }\KeywordTok{c}\NormalTok{( }\DecValTok{-400} \OperatorTok{*}\StringTok{ }\NormalTok{x[}\DecValTok{1}\NormalTok{] }\OperatorTok{*}\StringTok{ }\NormalTok{(x[}\DecValTok{2}\NormalTok{] }\OperatorTok{-}\StringTok{ }\NormalTok{x[}\DecValTok{1}\NormalTok{] }\OperatorTok{*}\StringTok{ }\NormalTok{x[}\DecValTok{1}\NormalTok{]) }\OperatorTok{-}\StringTok{ }\DecValTok{2} \OperatorTok{*}\StringTok{ }\NormalTok{(}\DecValTok{1} \OperatorTok{-}\StringTok{ }\NormalTok{x[}\DecValTok{1}\NormalTok{]),}
\DecValTok{200} \OperatorTok{*}\StringTok{ }\NormalTok{(x[}\DecValTok{2}\NormalTok{] }\OperatorTok{-}\StringTok{ }\NormalTok{x[}\DecValTok{1}\NormalTok{] }\OperatorTok{*}\StringTok{ }\NormalTok{x[}\DecValTok{1}\NormalTok{]) ) )}
\NormalTok{\}}
\end{Highlighting}
\end{Shaded}

We also need initial values

\begin{Shaded}
\begin{Highlighting}[]
\NormalTok{x0 <-}\StringTok{ }\KeywordTok{c}\NormalTok{( }\FloatTok{-1.2}\NormalTok{, }\DecValTok{1}\NormalTok{ )}
\end{Highlighting}
\end{Shaded}

Before we run the minimization procedure, we need to specify which
algorithm we will use. That can be done as follows:

\begin{Shaded}
\begin{Highlighting}[]
\NormalTok{opts <-}\StringTok{ }\KeywordTok{list}\NormalTok{(}\StringTok{"algorithm"}\NormalTok{=}\StringTok{"NLOPT_LD_LBFGS"}\NormalTok{,}
\StringTok{"xtol_rel"}\NormalTok{=}\FloatTok{1.0e-8}\NormalTok{)}
\end{Highlighting}
\end{Shaded}

Here, we will use the L-BFGS algorithm. Now we are ready to run the
optimization procedure.

\begin{Shaded}
\begin{Highlighting}[]
\CommentTok{# solve Rosenbrock Banana function}
\NormalTok{res <-}\StringTok{ }\KeywordTok{nloptr}\NormalTok{( }\DataTypeTok{x0=}\NormalTok{x0,}
\DataTypeTok{eval_f=}\NormalTok{eval_f,}
\DataTypeTok{eval_grad_f=}\NormalTok{eval_grad_f,}
\DataTypeTok{opts=}\NormalTok{opts)}
\end{Highlighting}
\end{Shaded}

See the result

\begin{Shaded}
\begin{Highlighting}[]
\KeywordTok{print}\NormalTok{(res)}
\end{Highlighting}
\end{Shaded}

\begin{verbatim}
## 
## Call:
## nloptr(x0 = x0, eval_f = eval_f, eval_grad_f = eval_grad_f, opts = opts)
## 
## 
## 
## Minimization using NLopt version 2.4.2 
## 
## NLopt solver status: 1 ( NLOPT_SUCCESS: Generic success return value. )
## 
## Number of Iterations....: 56 
## Termination conditions:  xtol_rel: 1e-08 
## Number of inequality constraints:  0 
## Number of equality constraints:    0 
## Optimal value of objective function:  7.35727226897802e-23 
## Optimal value of controls: 1 1
\end{verbatim}

The function is optimized at (1,1) which is the ground truth.

\section{Example 2: Minimizing with inequality constraint without
gradients}\label{example-2-minimizing-with-inequality-constraint-without-gradients}

The problem to minimize is

\[
\min_{x\in \mathbb{R}^n} \sqrt{x_2}
\] \[
\text{s.t.} x_2 \geq 0\] \[
x_2 \geq (a_1x_1 + b_1)^3\] \[
x_2 \geq  (a_2 x_1 + b_2)^3
\]

with \(a_1 = 2\), \(b_1 = 0\), \(a_2 = -1\), and \(b_2 = 1\). We
re-arrange the constraints to have the form \(g(x) \leq 0\):

\[
(a_1x_1 + b_1)^3 - x_2 \leq 0
\] \[
(a_2 x_1 + b_2)^3 - x+2 \leq 0
\]

First, define the objective function

\begin{Shaded}
\begin{Highlighting}[]
\CommentTok{# objective function}
\NormalTok{eval_f0 <-}\StringTok{ }\ControlFlowTok{function}\NormalTok{( x, a, b )\{}
\KeywordTok{return}\NormalTok{( }\KeywordTok{sqrt}\NormalTok{(x[}\DecValTok{2}\NormalTok{]) )}
\NormalTok{\}}
\end{Highlighting}
\end{Shaded}

and constraints are

\begin{Shaded}
\begin{Highlighting}[]
\CommentTok{# constraint function}
\NormalTok{eval_g0 <-}\StringTok{ }\ControlFlowTok{function}\NormalTok{( x, a, b ) \{}
\KeywordTok{return}\NormalTok{( (a}\OperatorTok{*}\NormalTok{x[}\DecValTok{1}\NormalTok{] }\OperatorTok{+}\StringTok{ }\NormalTok{b)}\OperatorTok{^}\DecValTok{3} \OperatorTok{-}\StringTok{ }\NormalTok{x[}\DecValTok{2}\NormalTok{] )}
\NormalTok{\}}
\end{Highlighting}
\end{Shaded}

Define parameters

\begin{Shaded}
\begin{Highlighting}[]
\CommentTok{# define parameters}
\NormalTok{a <-}\StringTok{ }\KeywordTok{c}\NormalTok{(}\DecValTok{2}\NormalTok{,}\OperatorTok{-}\DecValTok{1}\NormalTok{)}
\NormalTok{b <-}\StringTok{ }\KeywordTok{c}\NormalTok{(}\DecValTok{0}\NormalTok{, }\DecValTok{1}\NormalTok{)}
\end{Highlighting}
\end{Shaded}

Now solve using NLOPT\_LN\_COBYLA without gradient information

\begin{Shaded}
\begin{Highlighting}[]
\CommentTok{# Solve using NLOPT_LN_COBYLA without gradient information}
\NormalTok{res1 <-}\StringTok{ }\KeywordTok{nloptr}\NormalTok{( }\DataTypeTok{x0=}\KeywordTok{c}\NormalTok{(}\FloatTok{1.234}\NormalTok{,}\FloatTok{5.678}\NormalTok{),}
\DataTypeTok{eval_f=}\NormalTok{eval_f0,}
\DataTypeTok{lb =} \KeywordTok{c}\NormalTok{(}\OperatorTok{-}\OtherTok{Inf}\NormalTok{,}\DecValTok{0}\NormalTok{),}
\DataTypeTok{ub =} \KeywordTok{c}\NormalTok{(}\OtherTok{Inf}\NormalTok{,}\OtherTok{Inf}\NormalTok{),}
\DataTypeTok{eval_g_ineq =}\NormalTok{ eval_g0,}
\DataTypeTok{opts =} \KeywordTok{list}\NormalTok{(}\StringTok{"algorithm"}\NormalTok{=}\StringTok{"NLOPT_LN_COBYLA"}\NormalTok{,}
\StringTok{"xtol_rel"}\NormalTok{=}\FloatTok{1.0e-8}\NormalTok{),}
\DataTypeTok{a =}\NormalTok{ a,}
\DataTypeTok{b =}\NormalTok{ b )}
\KeywordTok{print}\NormalTok{( res1 )}
\end{Highlighting}
\end{Shaded}

\begin{verbatim}
## 
## Call:
## nloptr(x0 = c(1.234, 5.678), eval_f = eval_f0, lb = c(-Inf, 0), 
##     ub = c(Inf, Inf), eval_g_ineq = eval_g0, opts = list(algorithm = "NLOPT_LN_COBYLA", 
##         xtol_rel = 1e-08), a = a, b = b)
## 
## 
## Minimization using NLopt version 2.4.2 
## 
## NLopt solver status: 4 ( NLOPT_XTOL_REACHED: Optimization stopped because 
## xtol_rel or xtol_abs (above) was reached. )
## 
## Number of Iterations....: 50 
## Termination conditions:  xtol_rel: 1e-08 
## Number of inequality constraints:  2 
## Number of equality constraints:    0 
## Optimal value of objective function:  0.544331053951819 
## Optimal value of controls: 0.3333333 0.2962963
\end{verbatim}

\section{Example 3: Minimization with equality and inequality
constraints without
gradients}\label{example-3-minimization-with-equality-and-inequality-constraints-without-gradients}

We want to solve the following constraint optimization problem

\[
\min_{x} x_1 x_4(x_1 + x_2 + x_3) + x_3\]\[
\text{s.t.}\]\[
    x_1 x_2 x_3 x_4 \geq 25\]\[
    x_1^2 + x_2^2 + x_3^2 + x_4^2 = 40\]\[
    1 <= x_1,x_2,x_3,x_4 \leq 5
\]

Let's first solve it with gradients

\begin{Shaded}
\begin{Highlighting}[]
\NormalTok{eval_f <-}\StringTok{ }\ControlFlowTok{function}\NormalTok{( x ) \{}
        \KeywordTok{return}\NormalTok{( }\KeywordTok{list}\NormalTok{( }\StringTok{"objective"}\NormalTok{ =}\StringTok{ }\NormalTok{x[}\DecValTok{1}\NormalTok{]}\OperatorTok{*}\NormalTok{x[}\DecValTok{4}\NormalTok{]}\OperatorTok{*}\NormalTok{(x[}\DecValTok{1}\NormalTok{] }\OperatorTok{+}\StringTok{ }\NormalTok{x[}\DecValTok{2}\NormalTok{] }\OperatorTok{+}\StringTok{ }\NormalTok{x[}\DecValTok{3}\NormalTok{]) }\OperatorTok{+}\StringTok{ }\NormalTok{x[}\DecValTok{3}\NormalTok{],}
                      \StringTok{"gradient"}\NormalTok{ =}\StringTok{ }\KeywordTok{c}\NormalTok{( x[}\DecValTok{1}\NormalTok{] }\OperatorTok{*}\StringTok{ }\NormalTok{x[}\DecValTok{4}\NormalTok{] }\OperatorTok{+}\StringTok{ }\NormalTok{x[}\DecValTok{4}\NormalTok{] }\OperatorTok{*}\StringTok{ }\NormalTok{(x[}\DecValTok{1}\NormalTok{] }\OperatorTok{+}\StringTok{ }\NormalTok{x[}\DecValTok{2}\NormalTok{] }\OperatorTok{+}\StringTok{ }\NormalTok{x[}\DecValTok{3}\NormalTok{]),}
\NormalTok{                                      x[}\DecValTok{1}\NormalTok{] }\OperatorTok{*}\StringTok{ }\NormalTok{x[}\DecValTok{4}\NormalTok{],}
\NormalTok{                                      x[}\DecValTok{1}\NormalTok{] }\OperatorTok{*}\StringTok{ }\NormalTok{x[}\DecValTok{4}\NormalTok{] }\OperatorTok{+}\StringTok{ }\FloatTok{1.0}\NormalTok{,}
\NormalTok{                                      x[}\DecValTok{1}\NormalTok{] }\OperatorTok{*}\StringTok{ }\NormalTok{(x[}\DecValTok{1}\NormalTok{] }\OperatorTok{+}\StringTok{ }\NormalTok{x[}\DecValTok{2}\NormalTok{] }\OperatorTok{+}\StringTok{ }\NormalTok{x[}\DecValTok{3}\NormalTok{]) ) ) )}
\NormalTok{    \}}

    \CommentTok{# Inequality constraints.}
\NormalTok{    eval_g_ineq <-}\StringTok{ }\ControlFlowTok{function}\NormalTok{( x ) \{}
\NormalTok{        constr <-}\StringTok{ }\KeywordTok{c}\NormalTok{( }\DecValTok{25} \OperatorTok{-}\StringTok{ }\NormalTok{x[}\DecValTok{1}\NormalTok{] }\OperatorTok{*}\StringTok{ }\NormalTok{x[}\DecValTok{2}\NormalTok{] }\OperatorTok{*}\StringTok{ }\NormalTok{x[}\DecValTok{3}\NormalTok{] }\OperatorTok{*}\StringTok{ }\NormalTok{x[}\DecValTok{4}\NormalTok{] )}

\NormalTok{        grad   <-}\StringTok{ }\KeywordTok{c}\NormalTok{( }\OperatorTok{-}\NormalTok{x[}\DecValTok{2}\NormalTok{]}\OperatorTok{*}\NormalTok{x[}\DecValTok{3}\NormalTok{]}\OperatorTok{*}\NormalTok{x[}\DecValTok{4}\NormalTok{],}
                     \OperatorTok{-}\NormalTok{x[}\DecValTok{1}\NormalTok{]}\OperatorTok{*}\NormalTok{x[}\DecValTok{3}\NormalTok{]}\OperatorTok{*}\NormalTok{x[}\DecValTok{4}\NormalTok{],}
                     \OperatorTok{-}\NormalTok{x[}\DecValTok{1}\NormalTok{]}\OperatorTok{*}\NormalTok{x[}\DecValTok{2}\NormalTok{]}\OperatorTok{*}\NormalTok{x[}\DecValTok{4}\NormalTok{],}
                     \OperatorTok{-}\NormalTok{x[}\DecValTok{1}\NormalTok{]}\OperatorTok{*}\NormalTok{x[}\DecValTok{2}\NormalTok{]}\OperatorTok{*}\NormalTok{x[}\DecValTok{3}\NormalTok{] )}
        \KeywordTok{return}\NormalTok{( }\KeywordTok{list}\NormalTok{( }\StringTok{"constraints"}\NormalTok{=constr, }\StringTok{"jacobian"}\NormalTok{=grad ) )}
\NormalTok{    \}}

    \CommentTok{# Equality constraints.}
\NormalTok{    eval_g_eq <-}\StringTok{ }\ControlFlowTok{function}\NormalTok{( x ) \{}
\NormalTok{        constr <-}\StringTok{ }\KeywordTok{c}\NormalTok{( x[}\DecValTok{1}\NormalTok{]}\OperatorTok{^}\DecValTok{2} \OperatorTok{+}\StringTok{ }\NormalTok{x[}\DecValTok{2}\NormalTok{]}\OperatorTok{^}\DecValTok{2} \OperatorTok{+}\StringTok{ }\NormalTok{x[}\DecValTok{3}\NormalTok{]}\OperatorTok{^}\DecValTok{2} \OperatorTok{+}\StringTok{ }\NormalTok{x[}\DecValTok{4}\NormalTok{]}\OperatorTok{^}\DecValTok{2} \OperatorTok{-}\StringTok{ }\DecValTok{40}\NormalTok{ )}

\NormalTok{        grad   <-}\StringTok{ }\KeywordTok{c}\NormalTok{(  }\FloatTok{2.0}\OperatorTok{*}\NormalTok{x[}\DecValTok{1}\NormalTok{],}
                      \FloatTok{2.0}\OperatorTok{*}\NormalTok{x[}\DecValTok{2}\NormalTok{],}
                      \FloatTok{2.0}\OperatorTok{*}\NormalTok{x[}\DecValTok{3}\NormalTok{],}
                      \FloatTok{2.0}\OperatorTok{*}\NormalTok{x[}\DecValTok{4}\NormalTok{] )}
        \KeywordTok{return}\NormalTok{( }\KeywordTok{list}\NormalTok{( }\StringTok{"constraints"}\NormalTok{=constr, }\StringTok{"jacobian"}\NormalTok{=grad ) )}
\NormalTok{    \}}

    \CommentTok{# Initial values.}
\NormalTok{    x0 <-}\StringTok{ }\KeywordTok{c}\NormalTok{( }\DecValTok{1}\NormalTok{, }\DecValTok{5}\NormalTok{, }\DecValTok{5}\NormalTok{, }\DecValTok{1}\NormalTok{ )}

    \CommentTok{# Lower and upper bounds of control.}
\NormalTok{    lb <-}\StringTok{ }\KeywordTok{c}\NormalTok{( }\DecValTok{1}\NormalTok{, }\DecValTok{1}\NormalTok{, }\DecValTok{1}\NormalTok{, }\DecValTok{1}\NormalTok{ )}
\NormalTok{    ub <-}\StringTok{ }\KeywordTok{c}\NormalTok{( }\DecValTok{5}\NormalTok{, }\DecValTok{5}\NormalTok{, }\DecValTok{5}\NormalTok{, }\DecValTok{5}\NormalTok{ )}

    \CommentTok{# Optimal solution.}
\NormalTok{    solution.opt <-}\StringTok{ }\KeywordTok{c}\NormalTok{(}\FloatTok{1.00000000}\NormalTok{, }\FloatTok{4.74299963}\NormalTok{, }\FloatTok{3.82114998}\NormalTok{, }\FloatTok{1.37940829}\NormalTok{)}

    \CommentTok{# Set optimization options.}
\NormalTok{    local_opts <-}\StringTok{ }\KeywordTok{list}\NormalTok{( }\StringTok{"algorithm"}\NormalTok{ =}\StringTok{ "NLOPT_LD_MMA"}\NormalTok{,}
                        \StringTok{"xtol_rel"}\NormalTok{  =}\StringTok{ }\FloatTok{1.0e-7}\NormalTok{ )}
\NormalTok{    opts <-}\StringTok{ }\KeywordTok{list}\NormalTok{( }\StringTok{"algorithm"}\NormalTok{   =}\StringTok{ "NLOPT_LD_AUGLAG"}\NormalTok{,}
                  \StringTok{"xtol_rel"}\NormalTok{    =}\StringTok{ }\FloatTok{1.0e-7}\NormalTok{,}
                  \StringTok{"maxeval"}\NormalTok{     =}\StringTok{ }\DecValTok{1000}\NormalTok{,}
                  \StringTok{"local_opts"}\NormalTok{  =}\StringTok{ }\NormalTok{local_opts,}
                  \StringTok{"print_level"}\NormalTok{ =}\StringTok{ }\DecValTok{0}\NormalTok{ )}

    \CommentTok{# Do optimization.}
\NormalTok{    res <-}\StringTok{ }\KeywordTok{nloptr}\NormalTok{( }\DataTypeTok{x0          =}\NormalTok{ x0,}
                   \DataTypeTok{eval_f      =}\NormalTok{ eval_f,}
                   \DataTypeTok{lb          =}\NormalTok{ lb,}
                   \DataTypeTok{ub          =}\NormalTok{ ub,}
                   \DataTypeTok{eval_g_ineq =}\NormalTok{ eval_g_ineq,}
                   \DataTypeTok{eval_g_eq   =}\NormalTok{ eval_g_eq,}
                   \DataTypeTok{opts        =}\NormalTok{ opts )}

\KeywordTok{print}\NormalTok{(res)}
\end{Highlighting}
\end{Shaded}

\begin{verbatim}
## 
## Call:
## 
## nloptr(x0 = x0, eval_f = eval_f, lb = lb, ub = ub, eval_g_ineq = eval_g_ineq, 
##     eval_g_eq = eval_g_eq, opts = opts)
## 
## 
## Minimization using NLopt version 2.4.2 
## 
## NLopt solver status: 4 ( NLOPT_XTOL_REACHED: Optimization stopped because 
## xtol_rel or xtol_abs (above) was reached. )
## 
## Number of Iterations....: 537 
## Termination conditions:  xtol_rel: 1e-07 maxeval: 1000 
## Number of inequality constraints:  1 
## Number of equality constraints:    1 
## Optimal value of objective function:  17.014017291835 
## Optimal value of controls: 1 4.743002 3.821147 1.379409
\end{verbatim}

This can be solved differently without gradient as follows:

The objective function is defined as

\begin{Shaded}
\begin{Highlighting}[]
\NormalTok{eval_f <-}\StringTok{ }\ControlFlowTok{function}\NormalTok{(x)}
\NormalTok{\{}
    \KeywordTok{return}\NormalTok{ (x[}\DecValTok{1}\NormalTok{]}\OperatorTok{*}\NormalTok{x[}\DecValTok{4}\NormalTok{]}\OperatorTok{*}\NormalTok{(x[}\DecValTok{1}\NormalTok{] }\OperatorTok{+}\NormalTok{x[}\DecValTok{2}\NormalTok{] }\OperatorTok{+}\StringTok{ }\NormalTok{x[}\DecValTok{3}\NormalTok{] ) }\OperatorTok{+}\StringTok{ }\NormalTok{x[}\DecValTok{3}\NormalTok{] )}
\NormalTok{\}}
\end{Highlighting}
\end{Shaded}

Inequality constraint can be defined as:

\[
25 - x_1 x_2 x_3 x_4 \leq 0
\]

\begin{Shaded}
\begin{Highlighting}[]
\NormalTok{eval_g_ineq <-}\StringTok{ }\ControlFlowTok{function}\NormalTok{(x)}
\NormalTok{\{}
    \KeywordTok{return}\NormalTok{ (}\DecValTok{25} \OperatorTok{-}\StringTok{ }\NormalTok{x[}\DecValTok{1}\NormalTok{]}\OperatorTok{*}\NormalTok{x[}\DecValTok{2}\NormalTok{]}\OperatorTok{*}\NormalTok{x[}\DecValTok{3}\NormalTok{]}\OperatorTok{*}\NormalTok{x[}\DecValTok{4}\NormalTok{])}
\NormalTok{\}}
\end{Highlighting}
\end{Shaded}

Equality constraint can be defined as

\begin{Shaded}
\begin{Highlighting}[]
\NormalTok{eval_g_eq <-}\StringTok{ }\ControlFlowTok{function}\NormalTok{(x)}
\NormalTok{\{}
    \KeywordTok{return}\NormalTok{ ( x[}\DecValTok{1}\NormalTok{]}\OperatorTok{^}\DecValTok{2} \OperatorTok{+}\StringTok{ }\NormalTok{x[}\DecValTok{2}\NormalTok{]}\OperatorTok{^}\DecValTok{2} \OperatorTok{+}\StringTok{ }\NormalTok{x[}\DecValTok{3}\NormalTok{]}\OperatorTok{^}\DecValTok{2} \OperatorTok{+}\StringTok{ }\NormalTok{x[}\DecValTok{4}\NormalTok{]}\OperatorTok{^}\DecValTok{2} \OperatorTok{-}\StringTok{ }\DecValTok{40}\NormalTok{ )}
\NormalTok{\}}
\end{Highlighting}
\end{Shaded}

Let's specify upper and lower bounds

\begin{Shaded}
\begin{Highlighting}[]
\NormalTok{lb <-}\StringTok{ }\KeywordTok{c}\NormalTok{(}\DecValTok{1}\NormalTok{,}\DecValTok{1}\NormalTok{,}\DecValTok{1}\NormalTok{,}\DecValTok{1}\NormalTok{)}
\NormalTok{ub <-}\StringTok{ }\KeywordTok{c}\NormalTok{(}\DecValTok{5}\NormalTok{,}\DecValTok{5}\NormalTok{,}\DecValTok{5}\NormalTok{,}\DecValTok{5}\NormalTok{)}
\end{Highlighting}
\end{Shaded}

Now, define initial values

\begin{Shaded}
\begin{Highlighting}[]
\NormalTok{x0 <-}\StringTok{ }\KeywordTok{c}\NormalTok{(}\DecValTok{1}\NormalTok{,}\DecValTok{5}\NormalTok{,}\DecValTok{5}\NormalTok{,}\DecValTok{1}\NormalTok{)}
\end{Highlighting}
\end{Shaded}

Define options

\begin{Shaded}
\begin{Highlighting}[]
\CommentTok{# Set optimization options.}
\NormalTok{local_opts <-}\StringTok{ }\KeywordTok{list}\NormalTok{( }\StringTok{"algorithm"}\NormalTok{ =}\StringTok{ "NLOPT_LD_MMA"}\NormalTok{, }\StringTok{"xtol_rel"}\NormalTok{ =}\StringTok{ }\FloatTok{1.0e-15}\NormalTok{ )}
\NormalTok{opts <-}\StringTok{ }\KeywordTok{list}\NormalTok{( }\StringTok{"algorithm"}\NormalTok{=}\StringTok{ "NLOPT_GN_ISRES"}\NormalTok{,}
\StringTok{"xtol_rel"}\NormalTok{=}\StringTok{ }\FloatTok{1.0e-15}\NormalTok{,}
\StringTok{"maxeval"}\NormalTok{=}\StringTok{ }\DecValTok{160000}\NormalTok{,}
\StringTok{"local_opts"}\NormalTok{ =}\StringTok{ }\NormalTok{local_opts,}
\StringTok{"print_level"}\NormalTok{ =}\StringTok{ }\DecValTok{0}\NormalTok{ )}
\end{Highlighting}
\end{Shaded}

We use NL\_OPT\_LD\_MMA for local optimization and NL\_OPT\_GN\_ISRES
for overall optimization. You can set the tolerance to extremely low to
get the best result. The number of iterations is set using maxeval.
Setting tolerance to low or the number of iterations to very high may
result in the best approximation at the cost of increased computation
time.

Finally, optimize

\begin{Shaded}
\begin{Highlighting}[]
\NormalTok{res <-}\StringTok{ }\KeywordTok{nloptr}\NormalTok{ ( }\DataTypeTok{x0 =}\NormalTok{ x0,}
                \DataTypeTok{eval_f =}\NormalTok{ eval_f,}
                \DataTypeTok{lb =}\NormalTok{ lb,}
                \DataTypeTok{ub =}\NormalTok{ ub,}
                \DataTypeTok{eval_g_ineq =}\NormalTok{ eval_g_ineq,}
                \DataTypeTok{eval_g_eq =}\NormalTok{ eval_g_eq,}
                \DataTypeTok{opts =}\NormalTok{ opts}
\NormalTok{)}
\KeywordTok{print}\NormalTok{(res)}
\end{Highlighting}
\end{Shaded}

\begin{verbatim}
## 
## Call:
## 
## nloptr(x0 = x0, eval_f = eval_f, lb = lb, ub = ub, eval_g_ineq = eval_g_ineq, 
##     eval_g_eq = eval_g_eq, opts = opts)
## 
## 
## Minimization using NLopt version 2.4.2 
## 
## NLopt solver status: 5 ( NLOPT_MAXEVAL_REACHED: Optimization stopped because 
## maxeval (above) was reached. )
## 
## Number of Iterations....: 160000 
## Termination conditions:  xtol_rel: 1e-15 maxeval: 160000 
## Number of inequality constraints:  1 
## Number of equality constraints:    1 
## Current value of objective function:  17.1290602910051 
## Current value of controls: 1 4.45273 4.166179 1.347646
\end{verbatim}

\section{Example 4: Minimization with multiple inequality constraints
without
gradients}\label{example-4-minimization-with-multiple-inequality-constraints-without-gradients}

Our objective function in this case is

\[
\min_{x} x_1^2 + x_2^2
\]

such that

\[
1- x_1 + x_2 \leq 0
\] \[
1- x_1^2 + x_2^2 \leq 0
\] \[
9x_1^2 + x_2^2 \geq 9
\] \[
x_1^2 -x_2 >=0
\] \[
x_2^2 - x_1 \geq 0
\] with bounds on variables as

\[
-50 \leq x_1, x_2 \leq 50
\]

Let's write objective function first

\begin{Shaded}
\begin{Highlighting}[]
\NormalTok{eval_f <-}\StringTok{ }\ControlFlowTok{function}\NormalTok{(x)}
\NormalTok{\{}
    \KeywordTok{return}\NormalTok{ ( x[}\DecValTok{1}\NormalTok{]}\OperatorTok{^}\DecValTok{2} \OperatorTok{+}\StringTok{ }\NormalTok{x[}\DecValTok{2}\NormalTok{]}\OperatorTok{^}\DecValTok{2}\NormalTok{ )}
\NormalTok{\}}
\end{Highlighting}
\end{Shaded}

Inequality constraints can be written as

\begin{Shaded}
\begin{Highlighting}[]
\NormalTok{eval_g_ineq <-}\StringTok{ }\ControlFlowTok{function}\NormalTok{ (x) \{}
\NormalTok{    constr <-}\StringTok{ }\KeywordTok{c}\NormalTok{(}\DecValTok{1} \OperatorTok{-}\StringTok{ }\NormalTok{x[}\DecValTok{1}\NormalTok{] }\OperatorTok{-}\StringTok{ }\NormalTok{x[}\DecValTok{2}\NormalTok{], }
                \DecValTok{1} \OperatorTok{-}\StringTok{ }\NormalTok{x[}\DecValTok{1}\NormalTok{]}\OperatorTok{^}\DecValTok{2} \OperatorTok{-}\StringTok{ }\NormalTok{x[}\DecValTok{2}\NormalTok{]}\OperatorTok{^}\DecValTok{2}\NormalTok{,}
                \DecValTok{9} \OperatorTok{-}\StringTok{ }\DecValTok{9}\OperatorTok{*}\NormalTok{x[}\DecValTok{1}\NormalTok{]}\OperatorTok{^}\DecValTok{2} \OperatorTok{-}\StringTok{ }\NormalTok{x[}\DecValTok{2}\NormalTok{]}\OperatorTok{^}\DecValTok{2}\NormalTok{,}
\NormalTok{                x[}\DecValTok{2}\NormalTok{] }\OperatorTok{-}\StringTok{ }\NormalTok{x[}\DecValTok{1}\NormalTok{]}\OperatorTok{^}\DecValTok{2}\NormalTok{,}
\NormalTok{                x[}\DecValTok{1}\NormalTok{] }\OperatorTok{-}\StringTok{ }\NormalTok{x[}\DecValTok{2}\NormalTok{]}\OperatorTok{^}\DecValTok{2}\NormalTok{)}
    \KeywordTok{return}\NormalTok{ (constr)}
\NormalTok{\}}
\end{Highlighting}
\end{Shaded}

Lower and upper bounds are defined as

\begin{Shaded}
\begin{Highlighting}[]
\NormalTok{lb <-}\StringTok{ }\KeywordTok{c}\NormalTok{(}\OperatorTok{-}\DecValTok{50}\NormalTok{, }\DecValTok{-50}\NormalTok{)}
\NormalTok{ub <-}\StringTok{ }\KeywordTok{c}\NormalTok{(}\DecValTok{50}\NormalTok{, }\DecValTok{50}\NormalTok{)}
\end{Highlighting}
\end{Shaded}

Initial values are

\begin{Shaded}
\begin{Highlighting}[]
\NormalTok{x0 <-}\StringTok{ }\KeywordTok{c}\NormalTok{(}\DecValTok{3}\NormalTok{, }\DecValTok{1}\NormalTok{)}
\end{Highlighting}
\end{Shaded}

Finally, define options for \texttt{nloptr}

\begin{Shaded}
\begin{Highlighting}[]
\NormalTok{ opts <-}\StringTok{ }\KeywordTok{list}\NormalTok{( }\StringTok{"algorithm"}\NormalTok{            =}\StringTok{ "NLOPT_GN_ISRES"}\NormalTok{,}
                  \StringTok{"xtol_rel"}\NormalTok{             =}\StringTok{ }\FloatTok{1.0e-15}\NormalTok{,}
               \StringTok{"maxeval"}\NormalTok{=}\StringTok{ }\DecValTok{160000}\NormalTok{,}
                  \StringTok{"tol_constraints_ineq"}\NormalTok{ =}\StringTok{ }\KeywordTok{rep}\NormalTok{( }\FloatTok{1.0e-10}\NormalTok{, }\DecValTok{5}\NormalTok{ ))}
\end{Highlighting}
\end{Shaded}

Optimize

\begin{Shaded}
\begin{Highlighting}[]
\NormalTok{res <-}\StringTok{ }\KeywordTok{nloptr}\NormalTok{(}
        \DataTypeTok{x0          =}\NormalTok{ x0,}
        \DataTypeTok{eval_f      =}\NormalTok{ eval_f,}
        \DataTypeTok{lb          =}\NormalTok{ lb,}
        \DataTypeTok{ub          =}\NormalTok{ ub,}
        \DataTypeTok{eval_g_ineq =}\NormalTok{ eval_g_ineq,}
        \DataTypeTok{opts        =}\NormalTok{ opts )}
\KeywordTok{print}\NormalTok{(res)}
\end{Highlighting}
\end{Shaded}

\begin{verbatim}
## 
## Call:
## 
## nloptr(x0 = x0, eval_f = eval_f, lb = lb, ub = ub, eval_g_ineq = eval_g_ineq, 
##     opts = opts)
## 
## 
## Minimization using NLopt version 2.4.2 
## 
## NLopt solver status: 5 ( NLOPT_MAXEVAL_REACHED: Optimization stopped because 
## maxeval (above) was reached. )
## 
## Number of Iterations....: 160000 
## Termination conditions:  xtol_rel: 1e-15 maxeval: 160000 
## Number of inequality constraints:  5 
## Number of equality constraints:    0 
## Current value of objective function:  1.99999999972382 
## Current value of controls: 1 1
\end{verbatim}

While We didn't present the examples with multiple equality constraints,
they are very similar to Example 4. However, be sure to select the
optimization algorithm as NLOPT\_GN\_ISRES.

\renewcommand\refname{References}
\bibliography{biblio.bib}

\end{document}